
\magnification=1200
\vsize=20cm\hsize=13.5cm\hoffset=-0.1mm
\parindent=0mm \parskip=6pt plus 2pt minus 2pt
\abovedisplayskip=4pt plus 1pt minus 1pt
\belowdisplayskip=4pt plus 1pt minus 1pt

\footline{\hfill\rm  -- \folio\ --\hfill}

\ifx\mubyte\undefined 
\font\hugebf=cmbx10 at 14.4pt

\font\bigbf=cmbx10 at 12pt

\font\petcap=cmcsc10

\font\sevenbf=cmbx10 at 7pt
\else
\font\hugebf=ec-lmbx10 at 14.4pt
-lmbxi10 at 14.4pt
\font\bigbf=ec-lmbx10 at 12pt
-lmbxi10 at 12pt
-lmbxi10
\font\petcap=ec-lmcsc10
-lmr10 at 7pt
-lmri10 at 7pt
\font\sevenbf=ec-lmbx10 at 7pt
\fi

\font\tenCal=eusm10
\font\sevenCal=eusm7
\font\fiveCal=eusm5
\newfam\Calfam
  \textfont\Calfam=\tenCal
  \scriptfont\Calfam=\sevenCal
  \scriptscriptfont\Calfam=\fiveCal

\font\tenmsa=msam10
\font\sevenmsa=msam7
\font\fivemsa=msam5
\newfam\msafam
  \textfont\msafam=\tenmsa
  \scriptfont\msafam=\sevenmsa
  \scriptscriptfont\msafam=\fivemsa

\font\tenmsb=msbm10
\font\sevenmsb=msbm7
\font\fivemsb=msbm5
\newfam\msbfam
  \textfont\msbfam=\tenmsb
  \scriptfont\msbfam=\sevenmsb
  \scriptscriptfont\msbfam=\fivemsb
\def\Bbb{\fam\msbfam\tenmsb}


\def\bC{{\Bbb C}}

\def\bN{{\Bbb N}}

\def\hexnbr#1{\ifnum#1<10 \number#1\else
 \ifnum#1=10 A\else\ifnum#1=11 B\else\ifnum#1=12 C\else
 \ifnum#1=13 D\else\ifnum#1=14 E\else\ifnum#1=15 F\fi\fi\fi\fi\fi\fi\fi}
\def\msatype{\hexnbr\msafam}
\def\msbtype{\hexnbr\msbfam}
\mathchardef\leqslant="3\msatype36  \let\le\leqslant
\mathchardef\geqslant="3\msatype3E  \let\ge\geqslant
\mathchardef\compact="3\msatype62
\mathchardef\complement="3\msatype7B
\mathchardef\stimes="3\msatype02
\mathchardef\smallsetminus="2\msbtype72   
\mathchardef\squigrightarrow="3\msatype20
\mathchardef\subsetneq="3\msbtype20
\mathchardef\supsetneq="3\msbtype21

\def\section#1#2{{\baselineskip=15pt\parindent=6.5mm\bigbf
\item{#1.} #2\vskip2pt}}

\def\pointir{\kern2pt\raise-0.5pt\hbox{{\bigbf-}\kern-2pt{\bigbf-\kern2pt}}}

\def\\{\hfil\break}

\def\build#1|#2|#3|{\mathrel{\mathop{\null#1}\limits^{#2}_{#3}}}

\def\square{{\hfill \hbox{
\vrule height 1.453ex  width 0.093ex  depth 0ex
\vrule height 1.5ex  width 1.3ex  depth -1.407ex\kern-0.1ex
\vrule height 1.453ex  width 0.093ex  depth 0ex\kern-1.35ex
\vrule height 0.093ex  width 1.3ex  depth 0ex}}}
\def\qed{\phantom{$\quad$}\hfill$\square$\medskip}

\let\Degree=\degree \def\degree{\hbox{$\Degree$}}
\def\lguil{\hbox{%
\raise1pt\hbox{$\scriptscriptstyle\langle\!\langle\kern1pt$}}}
\def\rguil{\hbox{%
\raise1pt\hbox{$\kern1pt\scriptscriptstyle\rangle\!\rangle$}}}
\def\?{\hbox{$\,$}}

\let\leq=\leqslant
\let\geq=\geqslant

\def\Raise#1{\raise8pt\hbox{$#1$}}

\def\statement{\bf}
\def\frac#1#2{{#1\over #2}}

\def\st[#1]{\rlap{\kern-12.5mm[#1]}}


\footline{\hfill}
\headline{\ifnum \pageno=1 \hfil\else
\ifodd \pageno {\sevenbf \hfil
Jean-Pierre Demailly, Institut Fourier Grenoble\hfil\folio}%
\else {\sevenbf \folio\hfil Precise error estimate of the Brent-McMillan
algorithm for Euler's constant\hfil} \fi\fi}

\centerline{\hugebf Precise error estimate of the Brent-McMillan}
\bigskip
\centerline{\hugebf algorithm for the computation of Euler's constant}
\bigskip
\centerline{\bf Jean-Pierre Demailly}
\centerline{Universit\'e de Grenoble Alpes, Institut Fourier}
\vskip1.5cm

{\bf Abstract.} Brent and McMillan introduced in 1980 a new algorithm for the computation of Euler's constant $\gamma$, based on the use of the Bessel functions $I_0(x)$ and $K_0(x)$. It is the fastest known algorithm for the computation of $\gamma$. The time complexity can still be improved by evaluating a certain divergent asymptotic expansion up to its minimal term. Brent-McMillan conjectured in 1980 that the error is of the same magnitude as the last computed term, and Brent-Johansson partially proved it in 2015. They also gave some numerical evidence for a more precise estimate of the error term. We find here an explicit expression of that optimal estimate, along with a complete self-contained formal proof and an even more precise error bound.
\medskip
{\bf Key-words.} Euler's constant, Bessel functions, elliptic integral, integration by parts, asymptotic expansion, Euler-Maclaurin formula
\medskip
{\bf MSC 2010.} 11Y60, 33B15, 33C10, 33E05
\bigskip

\section{0}{Introduction and main results}

Let $H_n=1+{1\over 2}+\cdots+{1\over n}$ denote as usual the partial sums
of the harmonic series. The algorithm introduced by Brent-McMillan [BM80] 
for the computation of Euler's constant $\gamma=\lim_{n\to+\infty}(H_n-\log n)$ 
is based on certain identities satisfied by the Bessel functions
$I_\alpha(x)$ and $K_0(x)\,$:
$$
I_\alpha(x)
=\sum_{n=0}^{+\infty} {x^{\alpha+2n}\over n!\,\Gamma(\alpha+n+1)},\qquad
K_0(x)=-{\partial I_\alpha(x)\over\partial\alpha}_{\,|\alpha=0}.
\leqno(0.1)
$$
Experts will observe that $2x$ has been substituted to $x$ in the
conventional notation of Watson's treatise [Wat44]. As we will check
in \S$\,$1, these functions
satisfy the relations
$$
\leqalignno{
K_0(x) &= - (\log x + \gamma)I_0(x) + S_0(x)\qquad\hbox{where}&(0.2)\cr
I_0(x) &= \sum_{n=0}^{+\infty} {x^{2n}\over n!^2},\qquad
S_0(x) = \sum_{n=1}^{+\infty} H_n\;{x^{2n}\over n!^2}.&(0.3)\cr}
$$
As a consequence, Euler's constant can be written as
$$
\gamma={S_0(x)\over I_0(x)}-\log x - {K_0(x)\over I_0(x)},\leqno(0.4)
$$
and one can show easily that
$$
0 < {K_0(x)\over I_0(x)} <\pi\;e^{-4x}\qquad\hbox{for $x\ge 1$}.\leqno(0.5)
$$
In the simpler version $(BM)$ of the algorithm proposed by Brent-McMillan, 
the remainder term $\smash{{K_0(x)\over I_0(x)}}$ is neglected~; a
precision $10^{-d}$ is then achieved for $x \simeq {1\over 4}\,(d\,\log 10
+\log\pi)$, and the power series $I_0(x)$, $S_0(x)$ must be summed up 
to $n = \lceil a_1x\rceil$ approximately, where $a_p$ is the unique 
positive root of the equation
$$
a_p(\log a_p - 1) = p.
\leqno(0.6)
$$
The calculation of
$$
I_0(x)=1+{x^2\over 1^2}\bigg(1+{x^2\over 2^2}\bigg(\cdots
{x^2\over (n-1)^2}\bigg(1+{x^2\over n^2}\bigg(\cdots\bigg)\bigg)\cdots
\bigg)\bigg)
$$
requires $2$ arithmetic operations for each term, and that of
$$
S_0(x)\simeq H_{0,N} - {1\over I_0(x)}\;
{x^2\over 1^2}\bigg(H_{1,N}+{x^2\over 2^2}\bigg(\cdots
{x^2\over (n-1)^2}\bigg(H_{n-1,N}+{x^2\over n^2}\bigg(H_{n,N}
+\cdots\bigg)\bigg)\cdots\bigg)\bigg)
$$
requires 4 operations. The time complexity of the algorithm $(BM)$ is thus
$$
BM(d) = a_1 \times {1\over 4}\,d\;\log 10 \times 6 \times d \simeq 12.4\,d^2.
\leqno(0.7)
$$
However, as in Sweeney's more elementary method [Swe63], Brent and Mcmillan
observed that the remainder term $K_0(x)/I_0(x)$ can be evaluated by means of a
divergent asymptotic expansion
$$
I_0(x)K_0(x)\sim {1\over 4x}
\sum_{k\in\bN}{(2k)!^3\over k!^4\,(16x)^{2k}}.
\leqno(0.8)
$$
Their idea is to truncate the asymptotic expansion precisely at the
minimal term, which turns out to be obtained for $k=2x$ if $x$ is a
positive integer.
We will check, as was conjectured by Brent-McMillan [BM80] and partly
proven by Brent and Johansson [BJ15], that the corresponding
``truncation error'' is then of an order of magnitude comparable to the
minimal term $k=2x$, namely
$\smash{e^{-4x}\over 2\,\sqrt{2\pi}\;x^{3/2}}$ by Stirling's formula.
\medskip

{\statement Theorem.} {\it The truncation error
$$
\Delta(x):=
I_0(x)K_0(x)-{1\over 4x}\sum_{k=0}^{2x}{(2k)!^3\over k!^4\,(16x)^{2k}}
\leqno(0.9)
$$
admits when $x\to+\infty$ an equivalent
$$
\Delta(x) \sim -{5\,e^{-4x}\over 24\sqrt{2\pi}\,x^{3/2}},
\leqno(0.10)
$$
and more specifically
$$
\Delta(x)=-e^{-4x}\bigg(
{5\over 24\sqrt{2\pi}\,x^{3/2}}+\varepsilon(x)\bigg),
\qquad|\varepsilon(x)|<{0.863\over x^2}.
\leqno(0.11)
$$
The approximate value
$$
{K_0(x)\over I_0(x)}\simeq {1\over 4x\,I_0(x)^2}
\sum_{k=0}^{2x}{(2k)!^3\over k!^4\,(16x)^{2k}}\leqno(0.12)
$$
is thus affected by an error of magnitude
$$
{\Delta(x)\over I_0(x)^2}\sim 
\smash{-{5\sqrt{2\pi}\over 12\,x^{1{/}2}}\,e^{-8x}}.\leqno(0.13)
$$}%

The refined version $(BM')$ of the Brent-McMillan algorithm consists in
evaluating the remainder term ${K_0(x)\over I_0(x)}$ up to the accuracy
$e^{-8x}$ permitted by the approximation~(0.13). This implies to take
$x = {1\over 8}\,d\,\log 10$ and leads to a time complexity
$$
BM'(d) = \Big({3\over 4}a_3 + {1\over 2}\Big)\,\log 10\;d^2\simeq 9.7\,d^2,
\leqno(0.14)
$$
substantially better than (0.7). The proof of the above theorem
requires many calculations. The techniques developed here would
probably even yield an asymptotic development for~$\Delta(x)$, at
least for the first few terms, but the required calculations seem
very extensive. Hopefully, further asymptotic 
expansions of the error might be useful to investigate the arithmetic 
properties of $\gamma$, especially its rationality or irrationality.

The present paper is an extended version of an original text [Dem85]
written in June 1984 and published in ``Gazette des Math\'ematiciens'' in
1985. However, because of length constraints for such a mainstream 
publication, the main idea for obtaining the error estimate of the
Brent-McMillan algorithm had only been hinted, and most of the details
had been omitted. After more than 30 years passed, we take the
opportunity to make these details available and to improve the recent 
results of Brent-Johansson [BJ15].\medskip

\section{1}{Proof of the basic identities}

Relations (0.2) and (0.3) are obtained by using a derivation term by term
of the series defining $I_\alpha(x)$ in (0.1), along with the standard
formula ${\Gamma'(n+1)\over\gamma(n+1)}=H_n-\gamma$,
itself a consequence of the equalities
$$
{\Gamma'(x+1)\over\Gamma(x+1)}={1\over x}+
{\Gamma'(x)\over\Gamma(x)}\quad\hbox{and}\quad
\Gamma'(1)=-\gamma.
$$
Explicitly, we get
$$
{\partial I_\alpha(x)\over\partial\alpha}=\sum_{n=0}^{+\infty}
{\log x \cdot x^{\alpha+2n}\over n!\,\Gamma(\alpha+n+1)}-
{\Gamma'(\alpha+n+1)\, x^{\alpha+2n}\over n!\,\Gamma(\alpha+n+1)^2},
\leqno(1.1)
$$
hence (0.2) and (0.3). Now, the Hankel integral formula (see [Art31])
expresses the function $1/\Gamma$ as
$$
{1\over\Gamma(z)} = {1\over 2\pi i}\int_{(C)}\zeta^{-z}e^\zeta\;d\zeta
\leqno(1.2)
$$
where  $(C)$ is the open contour formed by a small circle 
$\zeta=\varepsilon e^{iu}$, $u\in[-\pi,\pi]$, concatenated with 
two half-lines $]-\infty,-\varepsilon]$ with respective arguments 
$-\pi$ and $+\pi$ and opposite orientation. This formula gives
$$
\leqalignno{
I_\alpha(x)&=\sum_{n=0}^{+\infty}{x^{\alpha+2n}\over n!}
{1\over 2\pi i}\int_{(C)}\zeta^{-\alpha-n-1}e^\zeta\;d\zeta=
{1\over 2\pi i}\int_{(C)}x^\alpha \zeta^{-\alpha-1}\exp(x^2/\zeta+\zeta)\;d\zeta\cr
&={1\over 2\pi i}\int_{(C)}\zeta^{-\alpha}\exp(x/\zeta+\zeta x)\;d\zeta\cr
&={1\over\pi}\int_0^\pi e^{2x\cos u}\cos(\alpha u)\;du
-{\sin\alpha\pi\over\pi}\int_0^{+\infty}e^{-2x\cosh v}\,e^{-\alpha v}\,dv.
&(1.3)\cr}
$$
The integral expressing $I_\alpha(x)$ in the second line above is obtained
by means of a change of variable $\zeta\mapsto \zeta x$ (recall that $x>0$)~;
the first integral of the third line comes from the modified contour
consisting of the circle $\{\zeta=e^{iu}\}$ of center $0$ and radius~$1$, and 
the last integral comes from the corresponding two half-lines
$t\in{}]-\infty,-1]$ written as $t=-e^{-v}$, $v\in{}]0,+\infty[\,$.
In~particular, the following integral expressions and equivalents of
$I_0(x)$,~$K_0(x)$ hold when $x\to +\infty$~:
$$
\leqalignno{
I_0(x)&={1\over\pi}\int_0^\pi e^{2x\cos u}\;du\kern27.5pt
\hbox{hence}~~
I_0(x) \mathop{\sim}\limits_{x\to+\infty}~{1\over\sqrt{4\pi x}}\,e^{2x},
&(1.4)\cr
K_0(x)&=\int_0^{+\infty}e^{-2x\cosh v}\,dv\qquad
\hbox{hence}~~
K_0(x) \mathop{\sim}\limits_{x\to+\infty}~\sqrt{\pi\over 4x}\,e^{-2x}.
&(1.5)\cr}
$$
Furthermore, one has $I_0(x)>{1\over\sqrt{4\pi x}}\,e^{2x}$ if $x\geq 1$ and
$K_0(x) <\sqrt{\pi\over 4x}\,e^{-2x}$ if $x>0$. These estimates can be
checked by means of changes of variables
$$
\eqalign{
I_0(x)&={e^{2x}\over 2\pi\sqrt{x}}\int_0^{4x} {e^{-t}\over\sqrt{t(1-t/4x)}}\;dt,
\qquad t=2x(1-\cos u),\cr
K_0(x)&={e^{-2x}\over 2\sqrt{x}}\int_0^{+\infty}{e^{-t}\over\sqrt{t(1+t/4x)}}\;dt,
\qquad t=2x(\cosh v-1),\cr}
$$
along with the observation that $\smash{\int_0^{+\infty}{1\over \sqrt{t}}\,
e^{-t}\,dt}=\Gamma({1\over 2})=\sqrt{\pi}$~; the lower bound for $I_0(x)$
is obtained by the convexity inequality
$\smash{{1\over \sqrt{1-t/4x}}}\geq 1+t/8x$ and an integration by
parts of the term $\sqrt{t}\;e^{-t}$, which give
$$
\eqalign{
\int_0^{4x} {e^{-t}\over\sqrt{t(1-t/4x)}}\;dt&\geq
\Gamma({\textstyle{1\over 2}})+{1\over 8x}\Gamma({\textstyle{3\over 2}})
-\int_{4x}^{+\infty}\Big({1\over\sqrt{t}}+{\sqrt{t}\over 8x}
\Big)\,e^{-t}\,dt\cr
&\geq \sqrt{\pi}+{\sqrt{\pi}\over 16x}
-e^{-4x}\Big({3\over 4\sqrt{x}}+{1\over 32x\,\sqrt{x}}\Big)
>\sqrt{\pi}\cr}
$$
for $x\geq 1$. Inequality (0.5) is then obtained by combining these bounds.
Our starting point to evaluate $K_0(x)$ more accurately is to use the integral
formulas (1.4), (1.5) to express $I_0(x)K_0(x)$ as a double integral 
$$
I_0(x)K_0(x)={1\over2\pi}
\int_{\{-\pi<u<\pi\,,\,v>0\}}\exp\big(2x(\cos u -\cosh v)\big)\,du\,dv.
\leqno(1.6)
$$
A change of variables
$$
r\;e^{i\theta}=\sin^2\Big({u+iv\over 2}\Big) = 
{1\over 2}\big(1-\cos(u+iv)\big)=
{1\over 2}\big(1-\cos u\cosh v+i\sin u\sinh v\big)
$$
gives
$$
\eqalign{
&r={1\over 2}(\cosh v-\cos u),\qquad|1-r\,e^{i\theta}|=
\Big|\cos\Big({u+iv\over 2}\Big)\Big|^2,\cr
&r\,dr\,d\theta=\Big|\sin\Big({u+iv\over 2}\Big)
\cos\Big({u+iv\over 2}\Big)\Big|^2\,du\,dv
=r\;|1-r\,e^{i\theta}|\,du\,dv,\cr}
$$
therefore
$$
I_0(x)K_0(x)
={1\over2\pi}\int_0^{+\infty}\exp(-4xr)\;dr\int_0^{2\pi}{d\theta\over
|1-r\,e^{i\theta}|}.\leqno(1.7)
$$
Let us denote by
$$
{\alpha\choose k}={\alpha(\alpha-1)\cdots(\alpha-k+1)\over k!},\qquad
\alpha\in\bC
$$
the (generalized) binomial coefficients.
For $z=r\,e^{i\theta}$ and $|z|=r<1$ the binomial identity $(1-z)^{-1/2}=
\sum_{k=0}^{+\infty}{-{1\over 2}\choose k}\,(-z)^k$
combined with the Parseval-Bessel formula yields the expansion$\phantom{\big|}$
$$
\varphi(r):={1\over2\pi}\int_0^{2\pi}{d\theta\over|1-r\,e^{i\theta}|}=
\sum_{k=0}^{+\infty}w_k\,r^{2k}\qquad 
\hbox{for $0\leq r<1$,}\leqno(1.8)
$$
where the coefficient
$$
w_k:={-1/2\choose k}^2=\bigg({1\cdot 3\cdot 5\cdots(2k-1)\over 
2\cdot 4\cdot 6\cdots 2k}\bigg)^2={(2k)!^2\over2^{4k}\,k!^4}.\leqno(1.9)
$$
is closely related to the Wallis integral
$\smash{W_p=\int_0^{\pi/2}\sin^p x\;dx}$. Indeed, the easily established
induction relation $W_p={p-1\over p}W_{p-2}$ implies
$$
W_{2k}={1\cdot 3\cdot 5\cdots (2k-1)\over 2\cdot 4\cdot 6\cdots 2k}\;
{\pi\over 2},\qquad
W_{2k+1}={2\cdot 4\cdot 6\cdots 2k\over 3\cdot 5\cdots (2k+1)},
$$
whence $w_k=({2\over \pi}W_{2k})^2$. The relations
$W_{2k}W_{2k-1}={\pi\over 4k}$, $W_{2k}W_{2k+1}={\pi\over 2(2k+1)}$ together
with the monotonicity of $(W_p)$ imply $\sqrt{\pi\over 2(2k+1)}<W_{2k}<
\sqrt{\pi\over 4k}$, therefore
$$
{2\over \pi(2k+1)}<w_k<{1\over \pi k}.\leqno(1.10)
$$
The main new ingredient of our analysis for estimating $I_0(x)K_0(x)$
is the following integral formula derived from (1.7), (1.8)$\,:$
$$
I_0(x)K_0(x)=\int_0^{+\infty}e^{-4xr}\,\varphi(r)\,dr
\leqno(1.11)
$$
where
$$
\leqalignno{
\varphi(r)&=\sum_{k=0}^{+\infty}w_k\,r^{2k}\kern85pt\hbox{for $r<1$,}
&(1.12)\cr
\varphi(r)&={1\over r}\varphi\bigg({1\over r}\bigg)=
\sum_{k=0}^{+\infty}w_k\,r^{-2k-1}\qquad\hbox{for $r>1$.}
&(1.13)\cr}
$$
(The last identity can be seen immediately by applying the change of
variable $\theta\mapsto-\theta$ in~(1.8)).
It is also easily checked using (1.10) that one has an equivalent
$$
\varphi(r)\sim\sum_{k=1}^{+\infty}{r^{2k}\over\pi k}=
{1\over \pi}\,\log{1\over 1-r^2}\qquad\hbox{when $r\to 1-0$},
$$
in particular the integral (1.11) converges near $r=1$
(later, we will need a more precise approximation, but more
sophisticated arguments are required for this).
By an integration term by term on $[0,+\infty[$ of the series 
defining $\varphi(r)$, and by ignoring the fact that the series diverges 
for $r\geq 1$, one formally obtains a divergent asymptotic expansion
$$
I_0(x)K_0(x)\sim \sum_{k\in\bN}w_k\,{(2k)!\over (4x)^{2k+1}} \sim
{1\over 4x}
\sum_{k\in\bN}{(2k)!^3\over k!^4\,(16x)^{2k}}.
\leqno(1.14)
$$
If $x$ is an integer, the general term of this expansion achieves its
minimum exactly for $k = 2x$, since the ratio of the $k$-th and 
$(k-1)$-st terms is
$$
{(2k(2k-1))^3\over k^4\,(16x)^2}=\bigg({k\over 2x}\bigg)^2
\bigg(1-{1\over 2k}\bigg)^3<1\quad\hbox{iff $k\le 2x$}.
$$
As already explained in the introduction, the idea is to truncate
the asymptotic expansion precisely at $k=2x$, and to estimate
the truncation error. This can be done by means of our explicit
integral formula (1.11).
\medskip

\section{2}{Expression of the error in terms of elliptic integrals}

By (1.7) and the definition of~$\Delta(x)$ we have
$$
\Delta(x)=\int_0^{+\infty}e^{-4xr}\,\delta(r)\,dr
\leqno(2.1)
$$
where
$$
\delta(r):=\varphi(r)-\sum_{k=0}^{2x}w_k\,r^{2k},\qquad
\hbox{so that}\quad\delta(r)=\sum_{k=2x+1}^{+\infty}w_k\,r^{2k}\quad
\hbox{for}\quad r<1.
\leqno(2.2)
$$
For $r<1$, let us also observe that $\varphi(r)$ coincides with the 
elliptic integral of the first kind
${2\over\pi}\int_0^{\pi/2}(1-r^2\sin^2\theta)^{-1/2}\,d\theta$, as
follows again from the binomial formula and the expression of~$W_{2k}$.
We need to calculate the precise asymptotic behavior of $\varphi(r)$
when~$r\to 1$. This can be obtained by means of a well known identity 
which we recall below. By putting $t^2=1-r^2$, the change of
variable $u=\tan\theta$ gives
$$
\leqalignno{
\varphi(r)&={2\over\pi}\int_0^{\pi/2}(1-r^2\cos^2\theta)^{-1/2}\,d\theta
={2\over\pi}\int_0^{+\infty}{du\over\sqrt{(1+u^2)(t^2+u^2)}}\,du\cr
&={4\over\pi}\int_0^1{dv\over \sqrt{(1+v^2)(1+t^2v^2)}}
+{2\over\pi}\int_t^1{dv\over \sqrt{(1+v^2)(t^2+v^2)}}&(2.3)\cr}
$$
where the last line is obtained by splitting the integral
$\int_0^{+\infty}\ldots\,du$ on the 3 intervals $[0,t]$, $[t,1]$,
$[1,+\infty[$, and by performing the respective changes of variable
$u=vt$, \hbox{$u=v$}, $u=1/v$ (the first and third pieces being then 
equal). Thanks to the binomial formula, the first integral of 
line (2.3) admits a development as a convergent series
$$
{4\over\pi}\int_0^1{dv\over \sqrt{(1+v^2)(1+t^2v^2)}}={4\over\pi}
\sum_{k=0}^{+\infty}
c'_kt^{2k},\qquad 
c'_k={-1/2\choose k}\int_0^1{v^{2k}\,dv\over\sqrt{1+v^2}}.
$$
The second integral can be expressed as the sum of a double series when
we simultaneously expand both square roots~:
$$
{2\over\pi}\int_t^1{dv\over v\sqrt{1+v^2}\,\sqrt{(1+t^2/v^2)}}
={2\over\pi}\int_t^1\sum_{k,\ell\geq 0}{-1/2\choose\ell}v^{2\ell}\,
{-1/2\choose k}(t^2/v^2)^k\;
{dv\over v}.
$$
The diagonal part $k=\ell$ yields a logarithmic term
$$
{2\over\pi}\sum_{k=0}^{+\infty}
{-1/2\choose k}^2\;t^{2k}\;
\log{1\over t}={1\over\pi}\,\varphi(t)\,\log{1\over t^2},
$$
and the other terms can be collected in the form of an absolutely convergent
double series
$$
{2\over\pi}\sum_{k\neq\ell\geq 0}{-1/2\choose k}{-1/2\choose\ell}t^{2k}
\Bigg[{v^{2\ell-2k}\over 2\ell-2k}\Bigg]_t^1
={2\over\pi}\sum_{k\neq\ell\geq 0}{-1/2\choose k}{-1/2\choose\ell}\;
{t^{2k}-t^{2\ell}\over 2(\ell-k)}.
$$
After grouping the various powers $t$, the summation reduces to a power
series ${4\over\pi}\sum c''_kt^{2k}$ of radius of convergence~$1$, where
(due to the symmetry in $k,\ell$)
$$
c''_k=\sum_{0\leq\ell<+\infty,\,\ell\neq k}~{1\over2(\ell-k)}
{-1/2\choose k}
{-1/2\choose\ell}.
$$
In fact, we see a priori from (1.10) that
$$
|c'_k|\leq {1\over\sqrt{\pi k}}\,{1\over 2k+1}=O(k^{-3/2}),
\kern90pt
$$
and
$$
|c''_k|\leq{1\over 2\sqrt{\pi k}}\Bigg({1\over k}+
\sum_{0<\ell\neq k}{1\over|\ell-k|\,\sqrt{\pi\ell}}\Bigg)
=O\bigg({\log k\over k}\bigg).
$$
In total, if we put $t^2=1-r^2$, the above relation implies
$$
\varphi(r)={1\over\pi}\bigg(\varphi(t)\,\log{1\over t^2}+
4\sum_{k=0}^{+\infty}c_k\,t^{2k}\bigg),\qquad
c_k=c'_k+c''_k,\leqno(2.4)
$$
and this identity will produce an arbitrarily precise expansion of
$\varphi(r)$ when $r\to 1$. In order to compute the coefficients, 
we observe that
$$
c_k=c'_k+c''_k={-1/2\choose k}\,\alpha_k
$$
with
$$
\eqalign{
\alpha_k&=\int_0^1{v^{2k}\,dv\over\sqrt{1{+}v^2}}
+\int_1^{+\infty}\Bigg({v^{2k}\over\sqrt{1{+}v^2}}-\sum_{\ell=0}^k
{-1/2\choose\ell}\,v^{2k-2\ell-1}\Bigg)\,dv
+\sum_{\ell=0}^{k-1}{1\over2(\ell{-}k)}{-1/2\choose\ell}.\cr}
$$
A direct calculation gives
$$
c_0=\alpha_0=\int_0^1{dv\over\sqrt{1{+}v^2}}+
\int_1^{+\infty}\bigg({1\over\sqrt{1{+}v^2}}-{1\over v}\bigg)dv=\log 2.
$$
Next, if we write
$$
{v^{2k}\over\sqrt{1+v^2}}=v^{2k-1}\cdot {v\over\sqrt{1+v^2}},\qquad
(\sqrt{1+v^2})'={v\over\sqrt{1+v^2}}
$$
and integrate by parts after factoring $v^{2k-1}$, we get
$$
\eqalign{
\alpha_k&=\sum_{\ell=0}^{k-1}{1\over2(\ell-k)}{-1/2\choose\ell}
+\Big[v^{2k-1}\,\sqrt{1+v^2}\,\Big]_0^1-
\int_0^1(2k-1)\,v^{2k-2}\,\sqrt{1+v^2}\,dv\cr
&\qquad{}+\Bigg[v^{2k-1}\Bigg(\sqrt{1+v^2}-\sum_{\ell=0}^k{-1/2\choose\ell}\,
{v^{1-2\ell}\over 1-2\ell}\Bigg)\Bigg]_1^{+\infty}\cr
&\qquad{}-\int_1^{+\infty}(2k-1)\,v^{2k-2}\Bigg(\sqrt{1+v^2}-\sum_{\ell=0}^k
{-1/2\choose\ell}\,{v^{1-2\ell}\over 1-2\ell}\Bigg)\,dv.\cr}
$$
This suggests to calculate $\alpha_k+(2k-1)\alpha_{k-1}$ and to use the
simplification
$$
v^{2k-2}\,\sqrt{1+v^2}-{v^{2k-2}\over\sqrt{1+v^2}}={v^{2k}\over\sqrt{1+v^2}}.
$$
We then infer
$$
\eqalign{
\alpha_k&+(2k-1)\alpha_{k-1}=-(2k-1)\alpha_k+\sum_{\ell=0}^k
{-1/2\choose\ell}\,{1\over 1-2\ell}\cr
&+\int_1^{+\infty}\kern-10pt(2k-1)\,v^{2k-2}\Bigg(\!\sum_{\ell=0}^k
{-1/2\choose\ell}\,
\bigg({v^{1-2\ell}\over 1-2\ell}-v^{1-2\ell}\bigg)-\sum_{\ell=0}^{k-1}
{-1/2\choose\ell}\,
v^{-1-2\ell}\!\Bigg)dv\cr
&+2k\sum_{\ell=0}^{k-1}{1\over2(\ell-k)}{-1/2\choose\ell}
+(2k-1)\sum_{\ell=0}^{k-2}{1\over 2(\ell-(k-1))}{-1/2\choose\ell}.\cr}
$$
A change of indices $\ell=\ell'-1$ in the sums corresponding to $k-1$
then eliminates almost all terms. There only remains the term $\ell=k$ 
in the first summation, whence the induction relation
$$
2k\,\alpha_k+(2k-1)\alpha_{k-1}=-{-1/2\choose k}\,{1\over 2k-1},\quad
\hbox{i.e.}\quad
{\alpha_k\over{-1/2\choose k}}
-{\alpha_{k-1}\over{-1/2\choose k-1}}=-{1\over 2k(2k-1)}.
$$
We get in this way
$$
{c_k\over{-1/2\choose k}^2}={\alpha_k\over{-1/2\choose k}}
={\alpha_0\over 1}-\sum_{\ell=1}^k{1\over 2\ell(2\ell-1)}
=\log 2-\sum_{\ell=1}^{2k}{(-1)^{\ell-1}\over \ell}
$$
and the explicit expression
$$
c_k=w_k\,\Bigg(\log 2-
\sum_{\ell=1}^{2k}{(-1)^{\ell-1}\over \ell}\Bigg).
\leqno(2.5)
$$
The remainder of the alternating series expressing $\log 2$ is bounded by half 
of last calculated term, namely $1/4k$, thus according to (1.10) we have
$0<c_k<{1\over \pi^2k^2}$ if $k\geq 1$, and the radius of convergence of
the series is~$1$. From (1.11) and (2.4) we infer
as $r\to 1-0$ the well known expansion of the elliptic integral
$$
\varphi(r)={1\over\pi}\Bigg(
\sum_{k=0}^{+\infty}w_kt^{2k}\log{1\over t^2}
+4\sum_{k=0}^{+\infty}c_kt^{2k}\Bigg),\qquad t^2=1-r^2,\leqno(2.6)
$$
with
$$
w_0=1,~~w_1={1\over 4},~~w_2={9\over 64},~~
c_0=\log 2,~~c_1={1\over 4}\bigg(\log 2-{1\over 2}\bigg),
~~c_2={9\over 64}\bigg(\log 2-{7\over 12}\bigg).
$$
Let us compute explicitly the first terms of the asymptotic expansion at $r=1$
by putting $r=1+h$, $h\to 0$. For $r=1+h<1$ ($h<0$) we have 
$t^2=1-r^2=-2h-h^2=2|h|(1+h/2)$, where
$$
\eqalign{
&\log{1\over t^2}=\log{1\over 2|h|(1+h/2)}=\log{1\over |h|}-\log 2-{1\over 2}h
+{1\over 8}h^2+O(h^2),\cr
&\sum_{k=0}^{+\infty}w_kt^{2k}=1+{1\over 4}(-2h-h^2)+{9\over 64}(2h)^2+O(h^3),\cr
&4\sum_{k=0}^{+\infty}c_kt^{2k}=4\log 2+
\bigg(\log 2 -{1\over 2}\bigg)(-2h-h^2)+
{9\over 16}\bigg(\log 2 -{7\over 12}\bigg)(2h)^2+O(h^3),\cr}
$$
and
$$
\eqalign{
\varphi(1+h)&={1\over\pi}\Bigg(
\bigg(1-{1\over 2}h+{5\over 16}h^2+O(h^3)
\bigg)\bigg(\log{1\over |h|}-\log 2-{1\over 2}h+
{1\over 8}h^2+O(h^3)\bigg)\cr
&\kern67pt{}+4\log 2-\big(2\log 2-1\big)h+
\bigg({5\over 4}\log 2-{13\over 16}\bigg)h^2+O(h^3)\Bigg).\cr}
$$
If terms are written by decreasing order of magnitude, we get
$$
\leqalignno{
\varphi(1+h)&={1\over\pi}\Bigg(
\log{1\over |h|}+3\log 2-{1\over 2}h\log{1\over |h|}-
\bigg({3\over 2}\log 2-{1\over 2}\bigg)h\cr
&\kern67pt{}+{5\over 16}h^2\log{1\over |h|}+
\bigg({15\over 16}\log 2-{7\over 16}\bigg)h^2
+O\bigg(h^3\log{1\over|h|}\bigg)\Bigg).&(2.7)\cr}
$$
For $r=1+h>1$, the identity $\varphi(r)={1\over r}\varphi({1\over r})$
gives in a similar way
$$
\varphi(r)=
{1\over 1+h}\Bigg({1\over\pi}\sum_{k=0}^{+\infty}w_kt^{2k}\log{1\over t^2}
+\sum_{k=0}^{+\infty}c_kt^{2k}\Bigg),\quad t^2=1-{1\over r^2}=2h-3h^2+O(h^3).
$$
After a few simplifications, one can see that the expansion (2.7) is still
valid for $h>0$. Passing to the limit $r\to 0$, $t\to 1-0$ in (2.6)
implies the relation $\sum_{k\geq 0}c_k={\pi\over 4}$. The following Lemma 
will be useful.
\medskip

{\statement Lemma A.} {\it For $h>0$, the difference
$$
\leqalignno{
\rho(h)&=\varphi(1+h)-{1\over\pi}\Bigg(\log{1\over h}+3\log 2
-{1\over 2}h\,\log{1\over h}-\bigg({3\over 2}\log 2-{1\over 2}\bigg)h\Bigg)
&(2.8)\cr
&=\varphi(1+h)-{1\over 2\pi}\bigg((h-2)\log{h\over 8}+h\bigg)&(2.9)\cr}
$$
admits the upper bound
$$
|\rho(h)|\leq h^2\bigg(2+\log\bigg(1+{1\over h}\bigg)\bigg).\leqno(2.10)
$$}%

{\it Proof.} A use of the Taylor-Lagrange formula gives
$(1+h)^{-1}=1-h+\theta_1h^2$, $t^2=1-{1\over r^2}=2h-3\theta_2h^2$,
with $\theta_i\in{}]0,1[$, and we also find $t^2\leq 2h$ and
$$
\eqalign{
\log{1\over t^2}&=\log {r^2\over (r-1)(r+1)}=\log{1\over h}+2\log(1+h)
-\log\bigg(1+{h\over 2}\bigg)-\log 2\cr
&=\log{1\over h}-\log 2+{3\over 2}h-{7\over 8}\theta_3h^2,\quad\theta_2\in
{}]0,1[,\cr}
$$
while the remainder terms $\sum_{k\ge 2}w_kt^{2k}$ and $\sum_{k\ge 2}c_kt^{2k}$
are bounded respectively by
$$
{w_2t^4\over 1-t^2}\leq 4w_2r^2h^2\leq{225\over 256}h^2\quad\hbox{and}\quad
{c_2t^4\over 1-t^2}\leq 4c_2r^2h^2<{1\over 10}h^2\quad\hbox{if}~~
h\leq{1\over 4},~r=1+h\leq{5\over 4}.
$$
For $h\leq{1\over 4}$ we thus get an equality
$$
\eqalign{
\varphi(1+h)&={1\over\pi}(1-h+\theta_1h^2)\times\Bigg(\cr
&\qquad{}
\bigg(1+{1\over 4}(2h-3\theta_2h^2)+{225\over 256}\theta_4h^2\bigg)
\bigg(\log{1\over |h|}-\log 2+{3\over 2}h-{7\over 8}\theta_3h^2\bigg)\cr
&\qquad{}+4\log 2+\bigg(\log 2-{1\over 2}\bigg)(2h-3\theta_2h^2)
+{4\over 10}\theta_5h^2\Bigg)\cr}
$$
with $\theta_i\in{}]0,1[\,$. In order to estimate $\rho(h)$,
we fully expand this expression and replace each term by an
upper bound of its absolute value. For $h\leq{1\over 4}$, 
this shows that $|\rho(h)|\leq h^2(0.885\,\log{1\over h}+2.11)$,
so that (2.10) is satisfied. For $h\geq{1\over 4}$, we write
$$
\rho'(h)=\varphi'(1+h)-{1\over2\pi}\bigg(\log{h\over 8}+2-{2\over h}\bigg),
\qquad
\varphi'(r)=-\sum_{k=0}^{+\infty}(2k+1)w_k\,r^{-2k-2},
$$
and by (1.10) we get
$$
\sum_{k=0}^{+\infty}{2\over\pi}r^{-2k-2}<-\varphi'(r)<{1\over r^2}+\sum_{k=1}^{+\infty}{3k\over\pi k}r^{-2k-2}
<\sum_{k=0}^{+\infty}r^{-2k-2}={1\over r^2-1},
$$
therefore
$$
\eqalign{
&{2\over\pi}\,{1\over h(h+2)}<-\varphi'(1+h)<{1\over h(h+2)},\cr
&{1\over2\pi}\bigg(\log{8\over h}-2+{2\over h}-{2\pi\over h(h+2)}\bigg)<
\rho'(h)<{1\over2\pi}\bigg(\log{8\over h}-2+{2\over h+2}\bigg).\cr}
$$
This implies
$$
\eqalign{
&-1.72<{1\over 2\pi}\bigg(\log 4{-}2{+}{1\over 4}{-}{32\pi\over 9}\bigg)<
\rho'(h)<{1\over 2\pi}\bigg(\log 32{-}2{+}{8\over 9}\bigg)<1.51
~~\hbox{on}~\bigg[{1\over 4},2\bigg],\cr
&-{1\over 2\pi}\bigg(\log{h\over 8}+2\bigg)<
\rho'(h)<{1\over 2\pi}\bigg(\log 4 -{3\over 2}\bigg)<0\qquad\hbox{
on~~$[2,+\infty[\,$},\cr}
$$
therefore $|\rho'(h)|\leq{1\over 2\pi}(h-1-\log 8+2)\leq{1\over 2\pi}h$
for $h\in[2,+\infty[$. Since $\rho(2)\simeq 0.00249<{1\over\pi}$, we see
that $|\rho(h)|\leq {1\over 4\pi}h^2$, and this shows that (2.10) still holds
on~$[2,+\infty[\,$. A numerical calculation of $\rho(h)$ at sufficiently
close points in the interval $[{1\over 4},2]$ finally yields (2.10) on
that interval.\qed
\medskip

Now we split the integral (2.1) on the intervals $[0,1]$ and
$[1,+\infty[\,$, starting with the integral of~$\varphi$ on the 
interval~$[1,+\infty[\,$. The change of variable $r=1+t/4x$ provides
$$
\int_1^{+\infty}e^{-4xr}\,\varphi(r)\,dr = {e^{-4x}\over 4x}
\int_0^{+\infty}e^{-t}\,\varphi\Big(1+{t\over 4x}\Big)\,dt,\leqno(2.11)
$$
and Lemma A~$(2.9)$ yields for this integral an approximation
$$
\eqalign{
{e^{-4x}\over 8\pi x}&
\int_0^{+\infty}e^{-t}\bigg(\Big({t\over 4x}-2\Big)
\log{t\over 32x}+{t\over 4x}\bigg)\,dt\cr
&={e^{-4x}\over 8\pi x}\bigg(
\log(32x)\Big(2-{1\over 4x}\Big)+2\gamma
+{1\over 4x}\int_0^{+\infty}e^{-t}(t\log t+t)\,dt\bigg)\cr
&={e^{-4x}\over 4\pi x}\bigg(\log x+\gamma+5\log 2-
{\log x\over 8x}-{\gamma+5\log 2-2\over 8x}\bigg),\cr}
$$
with an error bounded by
$$
\eqalign{
{e^{-4x}\over 4x}\int_0^{+\infty}e^{-t}\bigg({t\over 4x}\bigg)^2&\bigg(
2+\log\bigg(1+{4x\over t}\bigg)\bigg)dt\cr
&={e^{-4x}\over 4x}\bigg({1\over 4x^2}+{1\over 16x^2}
\int_0^{+\infty}t^2\,e^{-t}\,\log{t+4x\over t}\,dt\bigg).\cr}
$$
Writing
$$
0<\log{t+4x\over t}=\log{4x\over t} +
\log\bigg(1+{t\over 4x}\bigg)
\leq \log{4x\over t} + {t\over 4x},
$$
we further see that
$$
\eqalign{
\int_0^{+\infty}t^2\,e^{-t}\,\log{t+4x\over t}\,dt&\leq
\int_0^{+\infty}t^2\,e^{-t}\,\bigg(\log{4x\over t}+{t\over 4x}\bigg)\,dt
=2\log 4x+{3\over 2x}+2\gamma-3.\cr}
$$
We infer
$$
\int_1^{+\infty}e^{-4xr}\,\varphi(r)\,dr =
{e^{-4x}\over 4\pi x}\bigg(\log x+\gamma+5\log 2-
{\log x\over 8x}\bigg)+{e^{-4x}\over 4x}\,R_1(x),
\leqno(2.12)
$$
with 
$$
|R_1(x)|<{\gamma+5\log 2-2\over 8\pi x}+{1\over 4x^2}+
{2\log 4x+{3\over 2x}+2\over 16x^2}<{0.483\over x}\quad
\hbox{if $\bN\ni x\geq 1$},\leqno(2.13)
$$
thanks to a numerical evaluation of the sequence in a suitable range.\bigskip

\section{3}{Estimate of the truncated asymptotic expansion}

We now estimate the two integrals
$\int_0^1e^{-4xr}\,\sum\limits_{k\geq 2x+1}w_k\,r^{2k}\,dr$,
$\int_1^{+\infty}e^{-4xr}\,\sum\limits_{k\leq 2x}w_k\,r^{2k}\,dr$.
By means of iterated integrations by parts, we get
$$
\leqalignno{
\int_0^1e^{-4xr}\,r^{2k}\,dr&=e^{-4x}\sum_{\ell=1}^{+\infty}
{(4x)^{\ell-1}\over (2k+1)\cdots(2k+\ell)},&(3.1)\cr
\int_1^{+\infty}e^{-4xr}\,r^{2k}\,dr&=
{e^{-4x}\over 4x}\bigg(1+
\sum_{\ell=1}^{2k}{2k(2k-1)\cdots(2k-\ell+1)\over(4x)^\ell}\bigg).&(3.2)\cr}
$$
Combining the identities (2.1), (2.2), (2.12), (3.1), (3.2) we find
$$\quad
\Delta(x)={e^{-4x}\over 4x}\Bigg({1\over\pi}\Big(\log x+\gamma+5\log 2\Big)
-{\log x\over 8\pi x}-\sum_{k=0}^{2x}w_k+S(x)+R_1(x)+R_2(x)\Bigg)\leqno(3.3)
$$
with
$$
~~S(x)=\sum_{k=2x+1}^{+\infty}
\sum_{\ell=1}^{2x-1}{w_k\,(4x)^\ell\over (2k+1)\cdots(2k+\ell)}
-\sum_{k=1}^{2x}
\sum_{\ell=1}^{2x-1}w_k\,{2k(2k-1)\cdots(2k-\ell+1)
\over(4x)^\ell},\leqno(3.4)
$$
and
$$
~~R_2(x)=
\sum_{k=2x+1}^{+\infty}\sum_{\ell=2x}^{+\infty}
{w_k\,(4x)^\ell\over (2k+1)\cdots(2k+\ell)}
-\sum_{k=1}^{2x}\sum_{\ell=2x}^{+\infty}
w_k\,{2k(2k-1)\cdots(2k-\ell+1)\over(4x)^\ell}\kern-4pt
\leqno(3.5)
$$
(In the final summation, terms of index $\ell> 2k$ are zero).
Formula (3.3) leads us to study the asymptotic expansion
of $\sum_{k=0}^{2x}w_k$. This development is easy to establish from (2.6)
(one could even calculate it at an arbitrarily large order).
\medskip

{\statement Lemma B.} {\it One has
$$
\leqalignno{
\qquad~~&w_k={1\over\pi k}\bigg(1-{1\over 2(2k-1)}+\varepsilon_k\bigg)
\quad\hbox{where}~~
{1\over 12k(2k-1)}<\varepsilon_k<{5\over 16k(2k-1)},~~k\geq 1,
&(3.6)\cr
&\sum_{k=0}^{2x}w_k={1\over\pi}\Big(\log x+5\log 2+\gamma\Big)
+R_3(x),\qquad {1\over 4\pi x}<R_3(x)<{19\over 48\pi x}.
&(3.7)\cr}
$$}%

{\it Proof.} The lower bound (3.6) is a consequence of the Euler-Maclaurin's
formula [Eul15] applied to the function $f(x)=\log{2x-1\over 2x}$.
This yields
$$
{1\over 2}\,\log w_k=\sum_{i=1}^kf(i)=C+\int_1^kf(x)\,dx+{1\over 2}f(k)+
\sum_{j=1}^p{b_{2j}\over (2j)!}\,f^{(2j-1)}(k)+\tilde R_p
$$
where $C$ is a constant, and where the remainder term $\tilde R_p$
is the product of the next term by a factor $[0,1]$, namely
$$
{b_{2p+2}\over (2p+2)!}\,f^{(2p+1)}(k)={2^{2p+1}\,b_{2p+2}\over(2p+1)(2p+2)}
\bigg({1\over (2k-1)^{2p+1}}-{1\over (2k)^{2p+1}}\bigg).
$$
We have here
$$
\eqalign{
\int_1^kf(x)\,dx&={1\over 2}(2k-1)\log(2k-1)-k\log k-(k-1)\,\log 2\cr
&=\bigg(k-{1\over 2}\bigg)\log\bigg(1-{1\over 2k}\bigg)
-{1\over 2}\,\log k+{1\over 2}\,\log 2\cr}
$$
and the constant $C$ can be computed by the Wallis formula. Therefore,
with $b_2={1\over 6}$, we have
$$
\eqalign{
\log w_k&=\log{1\over\pi k} + 2k\,\log\bigg(1-{1\over 2k}\bigg)+1
+2\theta\, b_2\bigg({1\over (2k-1)}-{1\over 2k}\bigg)\cr
&\geq \log{1\over\pi k}-{1\over 4k}-\sum_{\ell=3}^{+\infty}
{1\over \ell (2k)^{\ell-1}}
> \log{1\over\pi k}-{1\over 4k}-{1\over 3}\,{1\over (2k)^2}\,
{1\over 1-{1\over 2k}}.\cr}
$$
The inequality $e^{-x}\geq 1-x$ then gives
$$
w_k>{1\over\pi k}\bigg(1-{1\over 4k}-{1\over 6k(2k-1)}\bigg)
={1\over\pi k}\bigg(1-{1\over 2(2k-1)}+{1\over 12k(2k-1)}\bigg)
$$
and the lower bound (3.6) follows for all $k\geq 1$. In the other direction,
we get
$$
\log w_k<\log{1\over\pi k} - {1\over 4k }- {1\over 12k^2} - {1\over 32k^3}
+ {1\over 6k(2k-1)}
=\log{1\over\pi k} - {1\over 4k} + {1\over 12k^2(2k-1)} - {1\over 32k^3}
$$
and the inequality $e^{-x}\leq 1-x+{1\over 2}x^2$ implies 
$$
\leqalignno{
w_k&<{1\over\pi k}\Bigg(1- \bigg({1\over 4k} - {1\over 12k^2(2k-1)} + {1\over 32k^3}\bigg) + {1\over 2}\bigg({1\over 4k}\bigg)^2\Bigg)\cr
&&\hbox{whence (by a difference of polynomials and a reduction to the same denominator)}\cr
w_k&<{1\over\pi k}\Bigg(1- {1\over 2(2k-1)}+{5\over 16k(2k-1)}\Bigg)\qquad
\hbox{if $k\geq 3$.}\cr}
$$
One can check that the final inequality still holds for $k=1,2$, and this
implies the estimate (3.6). On the other hand, formula (2.6) yields
$$
\eqalign{
w_0+\sum_{k=1}^{+\infty}\Big(w_k-{1\over \pi k}\Big)r^{2k}
&=\varphi(r)-{1\over\pi}\,\log{1\over 1-r^2}\cr
&={1\over\pi}\,\big(\varphi(t)-1\big)\,\log{1\over 1-r^2}
+{4\over\pi}\,\log2+\sum_{k\geq 1}
c_k\,t^{2k}\cr}
$$
with $t=\sqrt{1-r^2}$ and $\varphi(t)=1+O(1-r^2)$. By passing to the limit
when $r\to 1-0$ and $t\to 0$, we thus get
$$
w_0+\sum_{k=1}^{+\infty}\Big(w_k-{1\over \pi k}\Big)={4\over\pi}\,\log2.
$$
We infer
$$
w_0+\sum_{k=1}^{2x}\Big(w_k-{1\over \pi k}\Big)-{4\over\pi}\,\log2 
=\sum_{2x+1}^{+\infty}\Big({1\over \pi k}-w_k\Big)
$$
and the upper and lower bounds in (3.6) imply
$$
0<\sum_{2x+1}^{+\infty}\Big({1\over \pi k}-w_k\Big)\leq
\sum_{2x+1}^{+\infty}{1\over 2\pi\,k(2k-1)}
<\sum_{2x+1}^{+\infty}{1\over 4\pi}\,{1\over k(k-1)}
={1\over 8\pi x}.
$$
The Euler-Maclaurin estimate 
$$
\sum_{k=1}^{2x}{1\over k}=\log(2x)+\gamma+{1\over 4x}+{b_2\over 2(2x)^2}
- {b_4\over4(2x)^4}+\cdots\leqno(3.8)
$$
then finally yields (3.7).\qed

It remains to evaluate the sum $S(x)$. This is considerably more difficult,
as a consequence of a partial cancellation of positive and negative terms.
The approximation (3.6) obtained in Lemma~B implies
$$
S(x)={2\over\pi}\bigg(T(x)-{1\over 2}U(x)+{5\over 8}R_4(x)\bigg),\leqno(3.9)
$$
and if we agree as usual that the empty product 
$(2k{-}2)\cdots(2k{-}\ell{+}1)={1\over 2k-1}$ for $\ell=1$
is equal to~$1$, we get
$$
\leqalignno{
\qquad T(x)&=\sum_{\ell=1}^{2x-1}
\sum_{k=2x+1}^{+\infty}{(4x)^\ell\over 2k(2k+1)\cdots(2k+\ell)}
-\sum_{\ell=1}^{2x-1}
\sum_{k=1}^{2x}{(2k-1)\cdots(2k-\ell+1)\over(4x)^\ell},&(3.10)\cr
\qquad U(x)&=\sum_{\ell=1}^{2x-1}
\sum_{k=2x+1}^{+\infty}{(4x)^\ell\over (2k-1)\cdots(2k+\ell)}
-\sum_{\ell=1}^{2x-1}
\sum_{k=1}^{2x}{(2k-2)\cdots(2k-\ell+1)\over (4x)^\ell},&(3.11)\cr}
$$
where the new error term $R_4(x)$ admits the upper bound
$$
~~|R_4(x)|\leq
\sum_{\ell=1}^{2x-1}
\sum_{k=2x+1}^{+\infty}{(4x)^\ell/2k\over (2k-1)\cdots(2k+\ell)}
+\sum_{\ell=1}^{2x-1}
\sum_{k=1}^{2x}{(2k-2)\cdots(2k-\ell+1)\over 2k\,(4x)^\ell}.\leqno(3.12)
$$

\section{4}{Application of discrete integration by parts}

To evaluate the sums $T(x)$, $U(x)$ and $R_4(x)$, 
our method consists in performing first a summation over the index~$k$, and
for this, we use ``discrete integrations by parts''. Set
$$
u^{a,b}_k:={1\over (2k+a)(2k+a+1)\cdots(2k+b-1)},\qquad a\leq b\leqno(4.1)
$$
(agreeing that the denominator is~$1$ if $a=b$). Then
$$
\eqalign{
u^{a,b}_k-u^{a,b}_{k+1}&= {(2k+b)(2k+b+1)-(2k+a)(2k+a+1)
\over (2k+a)(2k+a+1)\cdots(2k+b+1)}\cr
&={(b-a)(4k+a+b+1)\over (2k+a)(2k+a+1)\cdots(2k+b+1)}.\cr}
$$
The inequalities $2(2k+a)\leq 4k+a+b+1\leq 2(2k+b+1)$ imply
$$
{1\over (2k+a+1)\cdots(2k+b+1)}\leq {u^{a,b}_k-u^{a,b}_{k+1}\over 2(b-a)}
\leq{1\over (2k+a)(2k+a+1)\cdots(2k+b)}
$$
with an upward error and a downward error both equal to
$$
{b-a+1\over 2}\,{1\over (2k+a)(2k+a+1)\cdots(2k+b+1)}.
$$
In particular, through a summation
$\sum_{k=2x+1}^{+\infty}{u^{a-1,b-1}_k-u^{a-1,b-1}_{k+1}\over 2(b-a)}$,
these inequalities imply
$$
\sum_{k=2x+1}^{+\infty}
{1\over (2k+a)\cdots(2k+b)}\leq 
{u^{a-1,b-1}_{2x+1}\over 2(b-a)}
={1\over 2(b-a)}\,{1\over (4x+a+1)\cdots(4x+b)},
$$
with an upward error equal to
$$
\eqalign{{b-a+1\over 2}
\sum_{k=2x+1}^{+\infty}
{1\over (2k+a-1)\cdots(2k+b)}
&\leq{1\over 4}\,{1\over (4x+a)\cdots(4x+b)}\cr}
$$
and an ``error on the error'' (again upwards) equal to
$$
{(b{-}a{+}1)(b{-}a{+}2)\over 4}
\sum_{k=2x+1}^{+\infty}
{1\over (2k+a-2)\cdots(2k+b)}
\leq{b{-}a{+}1\over 8}\,{1\over (4x+a-1)\cdots(4x+b)}.
$$
In other words, we find
$$
\leqalignno{
\sum_{k=2x+1}^{+\infty}{1\over (2k+a)\cdots(2k+b)}&=
{1\over 2(b-a)}{1\over (4x+a+1)\cdots(4x+b)}
-{1\over 4}{1\over (4x+a)\cdots(4x+b)}\cr
&\qquad\quad
{}+\theta\,{b-a+1\over 8}\,{1\over (4x+a-1)\cdots(4x+b)},\qquad
\theta\in[0,1].&(4.2^{a,b}_3)\cr}
$$
If necessary, one could of course push further this development to an 
arbitrary number of terms $p$ rather than $3$. We will denote the
corresponding expansion $(4.2^{a,b}_p)$, and will use it here in the
cases $p=2,3$. For the summations
$\smash{\sum_{k=1}^{2x}}\ldots~$, we similarly define
$$
v^{a,b}_k=(2k-a)(2k-a-1)\cdots(2k-b+1),\qquad a\leq b,\leqno(4.3)
$$
and obtain
$$
\eqalign{
v^{a,b}_k-v^{a,b}_{k-1}
&=(2k-a-2)\cdots(2k-b+1)\big((2k-a)(2k-a-1)-(2k-b)(2k-b-1)\big)\cr
&=(2k-a-2)\cdots(2k-b+1)\big((b-a)(4k-a-b-1)\big).\cr}
$$
For $a<b$, the inequalities $2(2k-b)\leq(4k-a-b-1)\leq 2(2k-a-1)$ imply
$$
(2k-a-2)\cdots(2k-b)\leq
{v^{a,b}_k-v^{a,b}_{k-1}\over 2(b-a)}\leq (2k-a-1)\cdots(2k-b+1)
$$
with an upward error and a downward error both equal to
$$
{1\over 2}(b-a-1)\,(2k-a-2)\cdots(2k-b+1).
$$
By considering the sum
$\sum_{k=1}^{2x}{v^{a,b}_k-v^{a,b}_{k-1}\over 2(b-a)}$, we obtain
$$
\sum_{k=1}^{2x}(2k-a-1)\cdots(2k-b+1)\geq {v^{a,b}_{2x}-
v^{a,b}_0\over 2(b-a)}
$$
with a downward error
$$
{b-a-1\over 2}
\sum_{k=1}^{2x}(2k-a-2)\cdots(2k-b+1)\leq{v^{a,b-1}_{2x}-v^{a,b-1}_0\over 4}
$$
and an upward error on the error equal to
$$
{(b-a-1)(b-a-2)\over 4}
\sum_{k=1}^{2x}(2k-a-2)\cdots(2k-b+2)
\leq{b-a-1\over 8}\,\Big(v^{a,b-2}_{2x}-v^{a,b-2}_0\Big),
$$
i.e.\ there exists $\theta\in[0,1]$ such that
$$
\leqalignno{
\sum_{k=1}^{2x}(2k-&a-1)\cdots(2k-b+1)\cr
\noalign{\vskip-8pt}\cr
&={1\over 2(b-a)}\Big(v_{2x}^{a,b}-v^{a,b}_0\Big)
+{1\over 4}\Big(v^{a,b-1}_{2x}-v^{a,b-1}_0\Big)
-\theta\,{b-a-1\over 8}\,\Big(v^{a,b-2}_{2x}-v^{a,b-2}_0\Big),\cr
&={1\over 2(b-a)}\,v_{2x}^{a,b}
+{1\over 4}\,v^{a,b-1}_{2x}-\theta\,{b-a-1\over 8}\,v^{a,b-2}_{2x}
+C^{a,b}_3,&(4.4^{a,b}_3)\cr}
$$
with
$$
|C^{a,b}_3|\leq {1\over 2(b-a)}|v^{a,b}_0|+{1\over 4}|v^{a,b-1}_0|+
{b-a-1\over 8}|v^{a,b-2}_0|,\leqno(4.5^{a,b}_3)
$$
especially $C^{a,b}_3=0$ if $a=0$. The simpler order 2 case (with an 
initial upward error) gives
$$
\leqalignno{
\sum_{k=1}^{2x}(2k-&a-2)\cdots(2k-b)
={1\over 2(b-a)}\Big(v_{2x}^{a,b}-v^{a,b}_0\Big)
-\theta\,{1\over 4}\Big(v^{a,b-1}_{2x}-v^{a,b-1}_0\Big)\cr
&={1\over 2(b-a)}\,(4x-a)\cdots(4x-b+1)
-\theta\,{1\over 4}\,(4x-a)\cdots(4x-b+2)+C^{a,b}_2.&(4.6^{a,b}_2)\cr}
$$
In the order 3 case, it will be convenient to use a further change
$$v^{a,b}_k-v^{a+1,b+1}_k=(2k-a-1)\cdots(2k-b+1)\Big((2k-a)-(2k-b)\Big)=
(b-a)v^{a+1,b}_k.
$$
If we apply this equality to the values $(a,b)$, $(a,b-1)$ and $k=2x$, we see
that the  $(4.4^{a,b}_3)$ development can be written in the equivalent form
$$
\leqalignno{
&\sum_{k=1}^{2x}(2k-a-1)\cdots(2k-b+1)-C^{a,b}_3\cr
\noalign{\vskip-8pt}\cr
&={1\over 2(b-a)}v_{2x}^{a+1,b+1}+{3\over 4}v^{a+1,b}_{2x}
+{b-a-1\over 8}\Big(2\,v^{a+1,b-1}_{2x}
-\theta\,v^{a,b-2}_{2x}\Big),\cr
&={1\over 2(b-a)}\,(4x-a-1)\cdots(4x-b)
+{3\over 4}\,(4x-a-1)\cdots(4x-b+1)\cr
&\kern1.1cm{}+{b-a-1\over 8} \Big(2(4x-a-1)\cdots(4x-b+2)
-\theta\,(4x-a)\cdots(4x-b+3)\Big)&(4.7^{\,a,b}_3)\cr}
$$
According to (3.10), $(4.2^{0,\ell}_3)$ and $(4.7^{\,0,\ell}_3)$, we
get
$$
T(x)=T'(x)-T''(x)+R_5(x)\leqno(4.8)
$$
with
$$
\leqalignno{
T'(x)&=
\sum_{\ell=1}^{2x-1}
{1\over 2\ell}\bigg({(4x)^\ell\over (4x+1)\cdots(4x+\ell)}-
{(4x-1)\cdots(4x-\ell)\over(4x)^\ell}\bigg),&(4.9)\cr
T''(x)&=
\sum_{\ell=1}^{2x-1}{1\over 4}\,{(4x)^\ell\over 4x(4x+1)\cdots(4x+\ell)}+
{3\over 4}\,{(4x-1)\cdots(4x-\ell+1)\over(4x)^\ell},&(4.10)\cr
\qquad |R_5(x)|&\leq{1\over 8}
\sum_{\ell=1}^{2x-1}\bigg({(\ell+1)\,(4x)^\ell\over (4x-1)4x\cdots(4x+\ell)}+
{2(\ell-1)\,(4x-1)\cdots(4x-\ell+2)\over(4x)^\ell}\bigg).&(4.11)\cr}
$$
The last term in the last line comes from formula $(4.7^{\,0,\ell}_3)$,
by observing that the inequalities $4x\leq 2(4x-\ell+2)$ $\ell\leq 2x-1$
imply
$$
4x(4x-1) \cdots(4x-\ell+3)\leq 2(4x-1)\cdots(4x-\ell+2).
$$
Similarly, thanks to (3.11), $(4.2^{-1,\ell}_2)$ and $(4.6^{0,\ell-1}_2)$,
we obtain the decomposition
$$
U(x)=U'(x)-U''(x)+R_6(x)\leqno(4.12)
$$
with
$$
\leqalignno{
\kern36pt U'(x)&=\sum_{\ell=1}^{2x-1}
{1\over 2(\ell+1)}\,{(4x)^\ell\over 4x\cdots(4x+\ell)}-
\sum_{\ell=2}^{2x-1}{1\over 2(\ell-1)}\,
{4x(4x-1)\cdots(4x-\ell+2)\over(4x)^\ell},&(4.13)\cr
U''(x)&={1\over 4x}\sum_{k=1}^{2x}{1\over 2k-1}\qquad\hbox{(negative
term $\ell=1$ appearing in $U(x)$)},&(4.14)\cr
\quad|R_6(x)|&\leq
{1\over 4}\sum_{\ell=1}^{2x-1}
{(4x)^\ell\over (4x-1)\cdots(4x+\ell)}+
{1\over 4}\sum_{\ell=2}^{2x-1}{4x(4x-1)\cdots(4x-\ell+3)\over(4x)^\ell}.
&(4.15)\cr}
$$
The remainder terms $R_2(x)$ $[\,$ resp.\ $R_4(x)\,$] can be bounded in 
the same way by means of $(4.2^{0,\ell}_2)$ and $(4.6^{-1,\ell-1}_2)$ $[\,$ resp.\
$(4.2^{-1,\ell}_2)$ and $(4.6^{0,\ell-2}_2)\,$]
and (1.10), (3.5), (3.12) lead to
$$
\leqalignno{
\kern20pt|R_2(x)|&\leq{2\over \pi}\Bigg(
\sum_{\ell=2x}^{+\infty}\sum_{k=2x+1}^{+\infty}
{(4x)^\ell\over (2k)\cdots(2k+\ell)}
+\sum_{\ell=2x}^{+\infty}\sum_{k=1}^{2x}
{(2k-1)\cdots(2k-\ell+1)\over(4x)^\ell}\Bigg)\cr
&\leq{2\over \pi}\sum_{\ell=2x}^{+\infty}{1\over 2\ell}\bigg(
{(4x)^\ell\over (4x+1)\cdots(4x+\ell)}+
{(4x+1)\cdots(4x-\ell+2)\over(4x)^\ell}\bigg),&(4.16)\cr
|R_4(x)|&\leq
\sum_{\ell=1}^{2x-1}
\sum_{k=2x+1}^{+\infty}{(4x)^{\ell-1}\over (2k-1)\cdots(2k+\ell)}
+\sum_{\ell=1}^{2x-1}
\sum_{k=1}^{2x}{(2k-2)\cdots(2k-\ell+2)\over (4x)^\ell}\cr
&\leq
\sum_{\ell=1}^{2x-1}{1\over 2(\ell+1)}\,{(4x)^{\ell-1}\over
4x\cdots(4x+\ell)}+\sum_{\ell=3}^{2x-1}{1\over 2(\ell-2)}
{4x(4x-1)\cdots(4x-\ell+3)\over(4x)^\ell}&(4.17)\cr
&+\sum_{k=1}^{2x}{1\over 2k(2k-1)}\,{1\over 4x}
+\sum_{k=1}^{2x}{1\over (2k-1)}\,{1\over (4x)^2}~~
\hbox{[terms $\ell=1,2$ in the summation]}.&(4.18)\cr}
$$
Finally, by (3.3), (3.7), (3.9) and (4.8), (4.12) we get
the decomposition
$$
\leqalignno{
\Delta(x)={e^{-4x}\over 4\pi x}\Bigg(
2\,T'(x)&-2\,T''(x)-U'(x)+U''(x)-{\log x\over 8x}\cr
&+\pi\Big(R_1(x)+R_2(x)-R_3(x)\Big)-{5\over 4}R_4(x)+2\,R_5(x)-R_6(x)\Bigg).
&(4.19)\cr}
$$
{\statement Lemma C.} {\it The following inequalities hold$\,:$
$$
\leqalignno{
&\log2-{1\over 8x}<
\sum_{k=1}^{2x}{1\over 2k(2k-1)}<\log2-{1\over 2(4x+1)},\kern90pt&(4.20)\cr
&\sum_{k=1}^{2x}{1\over 2k-1}
<{3\over 2}\log 2+{1\over 2}\Big(\log x+\gamma\Big)+{1\over 24x^2},&(4.21)\cr
&U''(x)={\log x\over 8x}+R_7(x),\qquad
0<R_7(x)<{1.37\over x}.&(4.22)\cr}
$$}

{\it Proof.} To check (4.20), we observe that the sum of the series is $\log 2$ 
and that the remainder of index $2x$ admits the upper bound
$$\eqalign{
{1\over 2(4x+1)}&=\sum_{k=2x+1}^{+\infty}{1\over 4}
\bigg({1\over k-1/2}-{1\over k+1/2}\bigg)\cr
&<\sum_{k=2x+1}^{+\infty}{1\over 2k(2k-1)}<\sum_{k=2x+1}{1\over 4}
\bigg({1\over k-1}-{1\over k}\bigg)={1\over 8x}.\cr}
$$
According to the Euler-Maclaurin expansion (3.8), we get on the one hand
$$
\eqalign{
\sum_{k=1}^{2x}{1\over 2k-1}&=\sum_{\ell=1}^{4x}{(-1)^{\ell-1}\over \ell}
+\sum_{\ell=1}^{2x}{1\over 2\ell}
=\sum_{k=1}^{2x}{1\over 2k(2k-1)}+{1\over 2}\sum_{k=1}^{2x}{1\over k}\cr
&<\log 2-{1\over 2(4x+1)}+
{1\over 2}\bigg(\log(2x)+\gamma+{1\over 4x}+{1\over 12(2x)^2}\bigg)\cr
&={3\over 2}\log 2+{1\over 2}\Big(\log x+\gamma\Big)+{1\over 8x(4x+1)}+
{1\over 96x^2},\cr}
$$
whence $(4.21)$, and on the other hand
$$
\eqalign{
\sum_{k=1}^{2x}{1\over 2k-1}
&>\log 2+{1\over 2}\bigg(\log(2x)+\gamma+{1\over 12(2x)^2}-
{1\over 120(2x)^4}\bigg)\cr
&>{3\over 2}\log 2+{1\over 2}\Big(\log x+\gamma\Big)+{1\over 96x^2}
-{1\over 1920x^4}.\cr}
$$
A straightforward numerical computation gives
${3\over 2}\log 2+{1\over 2}\gamma+{1\over 24}<1.37$, which
then implies~(4.22).\qed
\medskip

We will now check that all remainder terms $R_i(x)$ are of a lower order
of magnitude than the main terms, and in particular that they admit a
bound $O(1/x)$. The easier term to estimate is $R_6(x)$. One can indeed
use a very rough inequality
$$
|R_6(x)|\leq {1\over 4}\sum_{\ell=1}^{2x-1}{1\over 4x(4x-1)}+
{1\over 4}\sum_{\ell=2}^{2x-1}{1\over (4x)^2}\leq 
{1\over 4}\,{2x-1\over 4x(4x-1)}+
{1\over 4}\,{2x-2\over (4x)^2}<{1\over 16x}.\leqno(4.23)
$$
Consider now $R_4(x)$. We use Lemma~C to bound both summations 
appearing in (4.18), and get in this way
$$
[[(4.18)]]\leq {\log 2-{1\over 2(4x+1)}\over 4x}
+{{3\over 2}\log 2+{1\over 2}(\log x+\gamma)+{1\over 24x^2}\over(4x)^2}<{0.234\over x}
$$
(this is clear for $x$ large since ${1\over 4}\log 2 <0.234$ -- the precise
check uses a direct numerical calculation for smaller values of $x$). By even
more brutal estimates, we find
$$
\eqalign{
\sum_{\ell=1}^{2x-1}{1\over 2(\ell{+}1)}\,{(4x)^{\ell-1}\over
4x\cdots(4x+\ell)}&\leq
\sum_{\ell=1}^{2x-1}{1\over 2(\ell{+}1)}\,{1\over (4x)^2}
\leq{\log2x{+}\gamma{+}{1\over 4x}{+}{1\over 12(2x)^2}{-}1\over 32x^2}
<{0.025\over x},\cr
&\kern-122pt\sum_{\ell=3}^{2x-1}{1\over 2(\ell-2)}
{4x(4x-1)\cdots(4x-\ell+1)\over(4x)^{\ell+2}}\leq
\sum_{\ell=1}^{2x-3}{1\over\ell}\,{1\over 32x^2}
\leq{\log2x+\gamma\over 32x^2}<{0.040\over x}.\cr}
$$
This gives the final estimate
$$
|R_4(x)|\leq {0.299\over x}.\leqno(4.24)
$$

\section{5}{Further integral estimates}
In order to get an optimal bound of the other terms, and especially 
their differences, we are going to replace some summations by suitable
integrals. Before, we must estimate more precisely the partial products
$\prod(4x\pm j)$, and for this, we use the power series expansion of
their logarithms. For $t>0$, we have $t-{1\over 2}\,t^2<\log(1+t)<t$. 
By taking $t={j\over 4x}$, we~find
$$
-{\sum\limits_{1\leq j\leq\ell}~j\over 4x}<
\log{(4x)^\ell\over(4x+1)\cdots(4x+\ell)}=
\sum_{1\leq j\leq\ell}\log{1\over 1+{j\over 4x}}<
-{\sum\limits_{1\leq j\leq\ell}~j\over 4x}+
{\sum\limits_{1\leq j\leq\ell}~j^2\over 2(4x)^2}.
$$
Since $\sum_{1\leq j\leq\ell}~j={\ell(\ell+1)\over 2}$ and
$\sum_{1\leq j\leq\ell}~j^2={\ell(\ell+1)(2\ell+1)\over 6}$,
we get
$$
-{\ell(\ell{+}1)\over 8x}<\log{(4x)^\ell\over(4x{+}1)\cdots(4x{+}\ell)}
<-{\ell(\ell{+}1)\over 8x}+{\ell(\ell{+}1)(2\ell{+}1)\over 12\,(4x)^2},
$$
therefore
$$
\exp\bigg({1\over 32x}{-}{(\ell{+}1/2)^2\over 8x}\bigg)<
{(4x)^\ell\over(4x{+}1)\cdots(4x{+}\ell)}<\exp\bigg({1\over 32x}
{-}{(\ell{+}1/2)^2\over 8x}{+}{(\ell{+}1/2)^3\over 96x^2}\bigg).\kern-8pt
\leqno(5.1)
$$
For $\ell\leq 2x-1$ we have
$$
{(\ell+1/2)^2\over 8x}-{(\ell+1/2)^3\over 96x^2}
={(\ell+1/2)^2\over 8x}\bigg(1-{(\ell+1/2)\over 12x}\bigg)\geq
{5\over 6}\,{(\ell+1/2)^2\over 8x},
$$
hence (after performing a suitable numerical calculation)
$$
\leqalignno{
\qquad&{(4x)^\ell\over(4x+1)\cdots(4x+\ell)}<
\exp\bigg({1\over 32x}-{5\over 6}\,{(\ell+1/2)^2\over 8x}\bigg)
\kern53pt\hbox{for $\ell\leq 2x-1$,}
&(5.2)\cr
&{(4x)^\ell\over(4x+1)\cdots(4x+\ell)}<\exp\bigg({1\over 32x}-{5\over 6}
{(2x-1/2)^2\over 12x}\bigg)<{1.52\over x}\qquad\hbox{for $\ell\geq 2x-1$.}\cr}
$$
For $\ell\geq 2x$, each new factor is at most ${4x\over 4x+\ell}\leq
{2\over 3}$, thus
$$
\sum_{\ell=2x}^{+\infty}
{(4x)^\ell\over(4x+1)\cdots(4x+\ell)}<{1.52\over x}\sum_{p=1}^{+\infty}
\bigg({2\over 3}\bigg)^p<{3.04\over x}.
\leqno(5.3)
$$
On the other hand, the analogous inequality $-t-16\,t^2/26<\log(1-t)<-t$ applied with $t={j\over 4x}\leq 1/4$ implies
$$
-{\ell(\ell+1)\over 8x}-{16\,\ell(\ell+1)(2\ell+1)\over 6\cdot 26\,(4x)^2}
<\log{(4x-1)\cdots(4x-\ell)\over (4x)^\ell}<-{\ell(\ell+1)\over 8x}.
\leqno(5.4)
$$
As $\exp(1/4x)>1+1/4x$, we infer
$$
~~{(4x+1)\cdots(4x-\ell+2)\over(4x)^\ell}\leq
\bigg(1+{1\over 4x}\bigg)\,\exp\bigg(-{(\ell-1)(\ell-2)\over 8x}\bigg)
<\exp\bigg(-{\ell(\ell-3)\over 8x}\bigg),\kern-5pt\leqno(5.5)
$$
and the ratio of two consecutive upper bounds associated with
indices $\ell,\ell+1$ is less than $\exp(-(2\ell-2)/8x)\leq e^{-1/4}$ 
if $\ell=2x$ and less than $e^{-1/2}$ if $\ell\geq 2x+1$, thus
$$
\sum_{\ell=2x}^{+\infty}{(4x+1)\cdots(4x-\ell+2)\over(4x)^\ell}
\leq \exp\bigg({3\over 4}-{x\over 2}\bigg)\Bigg(1+e^{-1/4}\sum_{p=0}^{+\infty}
e^{-p/2}\Bigg)<{4.65\over x}.
$$
As $2\ell\geq 4x$, we deduce from (4.16) that
$$
|R_2(x)|\leq{2\over\pi}\,{1\over 4x}\,{7.69\over x}<{1.224\over x^2}\leqno(5.6)
$$
(but actually, one can see that $R_2(x)$ even decays exponentially).
By means of a standard integral-series comparison,
the inequalities (4.11), (5.2) and~(5.4) also provide
$$
\leqalignno{
|R_5(x)|&\leq{1\over 8}\sum_{\ell=1}^{2x-1}{\ell+1\over 4x(4x-1)}\,
\exp\bigg({1\over 32x}-{5\over 6}\,{(\ell+1/2)^2\over 8x}\bigg)
+2\,{\ell-1\over (4x)^2}\,\exp\bigg({3\ell\over 8x}-{\ell^2\over 8x}\bigg)\cr
&\leq{1\over 8(4x)(3x)}\Bigg(e^{1\over 32}\!\!
\int_0^{+\infty}\!\bigg(t\,{+}{3\over 2}\bigg)
\exp\bigg(\!-{5\over 6}{t^2\over 8x}\bigg)dt+{3\,e^{3\over 4}\over 2}\!
\int_0^{+\infty}t\,\exp\bigg(\!-{t^2\over 8x}\bigg)dt\!\Bigg)\cr
&={1\over 96x^2}\Bigg(e^{1\over 32}\bigg({24\over 5}x+{3\over 2}
\sqrt{48x\over 5}\,{1\over 2}\,\sqrt{\pi}\bigg)+6\,e^{3\over 4}\,x\Bigg)
<{0.229\over x}\qquad
\hbox{for $x\geq 1$}.&(5.7)\cr}
$$
It then follows from (3.9) and (5.1) that
$$
\eqalign{
T'(x)=\sum_{\ell=1}^{2x-1}{1\over 2\ell}\bigg(
&{(4x)^\ell\over(4x+1)\cdots(4x+\ell)}-
{(4x-1)\cdots(4x-\ell)\over(4x)^\ell}\bigg)\cr
&=\sum_{\ell=1}^{2x-1}{1\over2\ell}\,{(4x)^\ell\over(4x+1)\cdots(4x+\ell)}\,
\bigg(1-\prod_{j=1}^\ell\Big(1-{j\over 4x}\Big)\Big(1+{j\over 4x}\Big)\bigg)\cr
&\leq\sum_{\ell=1}^{2x-1}
\exp\bigg({1\over 32x}-{(\ell+1/2)^2\over 8x}+{(\ell+1/2)^3\over 96\,x^2}
\bigg)\,{(\ell+1)^2\over 96x^2}\;;\cr}
$$
to get this, we have used here the inequality $1-\prod(1-a_j)\leq\sum a_j$ 
with $a_j={j^2\over (4x)^2}<1$, and the identity
$\sum_{j\leq \ell}j^2={\ell(\ell+1)(2\ell+1)\over 6}$.
In the other direction, we have a lower bound 
$\prod(1-a_j)^{-1}-1\geq\sum a_j$, thus (5.3) implies
$$
\eqalign{
T'(x)=\sum_{\ell=1}^{2x-1}{1\over 2\ell}\,&
{(4x-1)\cdots(4x-\ell)\over(4x)^\ell}\Bigg(\prod_{j=1}^\ell
\bigg(1-\Big({j\over 4x}\Big)^2\bigg)^{-1}-1\Bigg)\cr
&\qquad{}\geq\sum_{\ell=1}^{2x-1}
\exp\bigg(-{\ell(\ell+1)\over 8x}-{(\ell+1/2)^3\over 78\,x^2}
\bigg)\,{(\ell+1)(2\ell+1)\over 12\,(4x)^2}\cr
&\qquad{}\geq\sum_{\ell=1}^{2x-1}
\exp\bigg(-{(\ell+1/2)^2\over 8x}-{(\ell+1/2)^3\over 78\,x^2}
\bigg)\,{(\ell+1)(\ell+1/2)\over 96x^2}\cr
&\qquad{}\geq\sum_{\ell=1}^{2x-1}
\exp\bigg(-{(\ell+1/2)^2\over 8x}\bigg)\bigg(1-{(\ell+1/2)^3\over 78\,x^2}
\bigg)\,{(\ell+1)(\ell+1/2)\over 96x^2}.\cr}
$$
We now evaluate these sums by comparing them to integrals. This gives
$$
T'(x)\leq e^{1\over 32x}\int_0^{2x}
\exp\bigg(-{t^2\over 8x}+{t^3\over 96x^2}\bigg)\,{(t+3/2)^2\over 96x^2}\,dt
$$
when we estimate the term of index $\ell$ by the corresponding integral
on the interval $[\ell-1/2,\ell+1/2]$. The change of variable
$$
u={t^2\over 8x}-{t^3\over 96x^2}={t^2\over 8x}\bigg(1-{t\over 12x}\bigg),
\qquad du={t\over 4x}\bigg(1-{t\over 8x}\bigg)dt
$$
implies $u\geq{5\over 48x}\,t^2$, hence
$t\leq\sqrt{48x\over 5}\,\sqrt{u}$. Moreover, a trivial convexity argument
yields $(1-{v\over p})^{-1}\leq 1+{1\over p-1}v$ if $v\leq 1\;$;
if we take $v={t\over 2x}$ and $p=6$ (resp.\ $p=3$), we find
$$
\eqalign{
&t=\sqrt{8xu}\,\bigg(1-{t\over 12x}\bigg)^{-1/2}
\leq\sqrt{8xu}\,\bigg(1+{t\over 20x}\bigg)
\leq\sqrt{8xu}\,\bigg(1+{\sqrt{3\over 125x}}\,\sqrt{u}\bigg),\cr
&dt={4x\over t}\bigg(1-{t\over 8x}\bigg)^{-1}du\leq
{4x\over t}\bigg(1+{t\over 6x}\bigg)du
\leq{4x\over \sqrt{8x u}\,
}\bigg(1+{2\over\sqrt{15x}}\,\sqrt{u}\bigg)du,\cr}
$$
therefore
$$
T'(x)\leq
{e^{1\over 32 x}\over 96x^2}\int_0^{+\infty}
e^{-u}\bigg({3\over 2}+\sqrt{8xu}\bigg(1+\sqrt{3\over 125x}\sqrt{u}\bigg)\bigg)^2
\bigg(1+{2\over\sqrt{15x}}\,\sqrt{u}\bigg)
{\sqrt{2x}\,du\over\sqrt{u}}.
$$
This integral can be evaluated evaluated explicitly, its dominant term being
equal to
$$
{e^{1\over 32 x}\over 96x^2}\int_0^{+\infty}
e^{-u}(\sqrt{8xu})^2{\sqrt{2x}\,du\over\sqrt{u}}
\sim{\sqrt {2}\over 12\sqrt{x}}\int_0^{+\infty}e^{-u}\,\sqrt{u}\,du
={\sqrt{2\pi}\over 24\,x^{1/2}}.
$$
Moreover, the factor $e^{1\over 32x}$ factor admits the (very rough!) 
upper bound $1+{1\over 31.5\,x}$, whence an error bounded by
$$
{\sqrt{2\pi}\over 24\,x^{1/2}}\cdot {1\over 31.5\,x}<{0.004\over x}.
$$
All other terms appearing in the integral involve terms
$O({1\over x})$ with coefficients which are  products of factors
$\Gamma(a)$, ${1\over 2}\leq a\le 2$, by coefficients whose sum
is bounded by
$$
{e^{1\over 32}\over 96}\bigg[
\bigg({3\over 2}+\sqrt{8}\bigg(1+\sqrt{3\over 125}\,\bigg)\bigg)^2
\bigg(1+{2\over\sqrt{15}}\bigg)\sqrt{2}-8\sqrt{2}\,\bigg]
<0.4021.
$$
As $\Gamma(a)\leq\sqrt{\pi}$, we obtain
$$
T'(x)<{\sqrt{2\pi}\over 24\,x^{1/2}}+{0.717\over x}.
$$
Similarly, one can obtain the following lower bound for $T'(x)\,$:
$$
\eqalign{
T'(x)&\geq\sum_{\ell=1}^{2x-1}
\exp\bigg(-{(\ell+1/2)^2\over 8x}\bigg)\,
\bigg(1-{(\ell+1/2)^3\over 78\,x^2}\bigg){(\ell+1)(\ell+1/2)\over 96x^2}\cr
&\geq\int_{3/2}^{2x+1/2}
\exp\bigg(-{t^2\over 8x}\bigg)\bigg(1-{t^3\over 78\,x^2}
\bigg)\,{(t-1)(t-1/2)\over 96x^2}\,dt\cr
&\geq\int_{2}^{2x}
\exp\bigg(-{t^2\over 8x}\bigg)\bigg(1-{t^3\over 78\,x^2}
\bigg)\,{t^2-3t/2\over 96x^2}\,dt\cr
&=\int_{1/2x}^{x/2}
e^{-u}\bigg(1-{8\sqrt{8}\,u^{3/2}\over 78\,x^{1/2}}
\bigg)\,{8xu-3\sqrt{8}\,x^{1/2}u^{1/2}/2\over 96x^2}\,{
\sqrt{8}\,x^{1/2}\,du\over 2\,u^{1/2}}\cr
&\geq\int_{1/2x}^{x/2}
e^{-u}\bigg(1-{8\sqrt{8}\,u^{3/2}\over 78\,x^{1/2}}
\bigg)\,{\sqrt{8}\,u-3\,x^{-1/2}u^{1/2}/2\over 24\,x^{1/2}}\,
{du\over u^{1/2}}\cr
&\geq\int_{1/2x}^{x/2}
e^{-u}\bigg({\sqrt{2}\,u^{1/2}\over 12\,x^{1/2}}-
{8\,u^2\over 3\cdot 78\,x}-{1\over 16x}\bigg)\,du\cr
&\geq\int_{0}^{+\infty}
e^{-u}\bigg({\sqrt{2}\,u^{1/2}\over 12\,x^{1/2}}-
{4\,u^2\over 117\,x}-{1\over 16x}\bigg)\,du-
\int_{\complement}
e^{-u}{\sqrt{2}\,u^{1/2}\over 12\,x^{1/2}}\,du.\cr}
$$
The integral $\int_{\complement}...$ on the ``missing intervals'' is bounded
on $[0,1/2x]$ by
$$
\int_{0}^{1/2x}{\sqrt{2}\,u^{1/2}\over 12\,x^{1/2}}\,du
={1\over 36\,x^2},
$$
whilst the integral on $[A,+\infty[{}=[x/2,+\infty[$ satisfies
$$
\int_A^{+\infty}u^\alpha\,e^{-u}\,du=A^\alpha\,e^{-A}+
\int_A^{+\infty}\alpha\,u^{\alpha-1}\,e^{-u}\,du
\leq e^{-A}(A^\alpha+\alpha A^{\alpha-1}),\quad\alpha\in{}]0,1].
$$
This provides an estimate
$$
\int_{{x\over 2}}^{+\infty}
e^{-u}{\sqrt{2}\,u^{1/2}\over 12\,x^{1/2}}\,du
\leq\exp\bigg(-{x\over 2}\bigg)\bigg({1\over 12}+
{1\over 12x}\bigg)
\leq{{1\over 6}\,e^{-1/2}\over x}.
$$
Therefore, we obtain the explicit lower bound
$$
T'(x)>{\sqrt{2\pi}\over 24\,x^{1/2}}-\bigg({8\over 117}+{1\over 16}+{1\over 36}
+{1\over 6}e^{-1/2}\bigg){1\over x}>{\sqrt{2\pi}\over 24\,x^{1/2}}-
{0.260\over x}.
$$
In the same manner, but now without any compensation of terms and with
much simpler calculations, the estimates (4.11), (5.1), (5.3) provide
an upper bound
$$
T''(x)\leq{1\over 4x}\sum_{\ell=1}^{2x-1}
{1\over 4}\exp\bigg({32\over x}-{(\ell{+}1/2)^2\over 8x}+{(\ell{+}1/2)^3\over 
96x^2}\bigg)+{3\over 4}\exp\bigg({32\over x}-{(\ell{-}1/2)^2\over 8x}\bigg).
$$
By using integral estimates very similar to those already used, this gives
$$
\eqalign{
T''(x)&\leq{e^{32\over x}\over 4x}\Bigg(
{1\over 4}\int_0^{2x}\exp\bigg(-{t^2\over 8x}+{t^3\over 96x^2}\bigg)dt
+{3\over 4}\int_0^{2x}
\exp\bigg(-{t^2\over 8x}\bigg)dt\Bigg)+{3\over 16x}\cr
&\leq{e^{32\over x}\over 4x}\Bigg(
{1\over 4}\int_0^{+\infty}e^{-u}\bigg(1+{2\over\sqrt{15x}}\,\sqrt{u}\bigg)
{\sqrt{2x}\,du\over \sqrt{u}}
+{3\over 4}\int_0^{+\infty}e^{-u}\,{\sqrt{2x}\,du\over \sqrt{u}}
\Bigg)+{3\over 16x}\cr
&\leq{e^{32\over x}\over 4x}\int_0^{+\infty}e^{-u}\,{\sqrt{2x}\,du\over \sqrt{u}}
+{e^{32\over x}\over 4x}\,{1\over\sqrt{30}}+{3\over 16x}
<{\sqrt{2\pi}\over 4x^{1/2}}+{0.255\over x},\cr}
$$
and we get likewise a lower bound
$$
\eqalign{
T''(x)&\geq{1\over 4x}\sum_{\ell=1}^{2x-1}
{1\over 4}\exp\bigg(-{(\ell+1/2)^2\over 8x}\bigg)+
{3\over 4}\exp\bigg(-{(\ell-1/2)^2\over 8x}
-{(\ell-1/2)^3\over 78x^2}
\bigg)\cr
&\geq{1\over 4x}\Bigg({1\over 4}\int_{3/2}^{2x+1/2}
\exp\bigg(-{t^2\over 8x}\bigg)dt
+{3\over 4}\int_{1/2}^{2x-1/2}
\exp\bigg(-{t^2\over 8x}\bigg)\bigg(1-{t^3\over 78x^2}\bigg)dt\Bigg)\cr
&\geq{1\over 4x}\Bigg({1\over 4}\int_0^{2x}
\exp\bigg(-{t^2\over 8x}\bigg)dt
+{3\over 4}\int_{0}^{2x}
\exp\bigg(-{t^2\over 8x}\bigg)\bigg(1-{t^3\over 78x^2}\bigg)dt-{9\over 8}\Bigg)\cr
&\geq{1\over 4x}\Bigg(\int_{0}^{2x}
\exp\bigg(-{t^2\over 8x}\bigg)dt
-{3\over 4}\int_0^{+\infty}
\exp\bigg(-{t^2\over 8x}\bigg)\,{t^3\over 78x^2}\,dt-{9\over 8}\Bigg)\cr
&={1\over 4x}\Bigg(\int_0^{x/2}e^{-u}\,{\sqrt{2x}\,du\over\sqrt{u}}
-{1\over 104}\int_0^{+\infty}e^{-u}{u\,du\over 2}-{9\over 8}\Bigg)\cr
&\geq{1\over 4x}\Bigg(\int_0^{+\infty}e^{-u}\,{\sqrt{2x}\,du\over\sqrt{u}}
-{235\over 208}-2\,e^{-x/2}\Bigg)>{\sqrt{2\pi}\over 4x^{1/2}}-{0.586\over x}.
\cr}
$$
All this finally yields the estimate
$$
T'(x)-T''(x)=-{5\over 24}\,{\sqrt{2\pi}\over x^{1/2}}+R_8(x),\qquad
-{0.515\over x}<R_8(x)<{1.303\over x}.
\leqno(5.8)
$$
There only remains to evaluate $U'(x)$. According to (4.13), a change of 
variable $\ell=\ell'+1$ followed by a decomposition $4x=(4x-\ell)+\ell$ 
allows us to transform the second summation appearing in $U'(x)$ as
$$
\eqalign{
U'(x)&=\sum_{\ell=1}^{2x-1}
{1\over 2(\ell+1)}\,{(4x)^\ell\over 4x\cdots(4x+\ell)}-
\sum_{\ell=1}^{2x-2}{1\over 2\ell}\,
{4x(4x-1)\cdots(4x-\ell+1)\over(4x)^{\ell+1}}\cr
&=\sum_{\ell=1}^{2x-1}
{1\over 2(\ell+1)}\,{(4x)^{\ell-1}\over (4x+1)\cdots(4x+\ell)}-
\sum_{\ell=1}^{2x-2}{1\over 2\ell}\,
{(4x-1)\cdots(4x-\ell+1)(4x-\ell)\over(4x)^{\ell+1}}\cr
&\kern160pt{}-\sum_{\ell=1}^{2x-2}{1\over 2}\,
{(4x-1)\cdots(4x-\ell+1)\over(4x)^{\ell+1}}.\cr}
$$
Writing ${1\over \ell+1}={1\over \ell}-{1\over \ell(\ell+1)}$, one obtains
$$
U'(x)={1\over 4x}\,T'(x)-R_9(x)
$$
with
$$
\eqalign{
R_9(x)&=\sum_{\ell=1}^{2x-1}
{1\over 2\ell(\ell+1)}\,{(4x)^{\ell-1}\over (4x+1)\cdots(4x+\ell)}
+\sum_{\ell=1}^{2x-2}{1\over 2}\,
{(4x-1)\cdots(4x-\ell+1)\over(4x)^{\ell+1}}\cr
&\kern100pt{}-\bigg({1\over 2\ell}{(4x-1)\cdots(4x-\ell)
\over(4x)^{\ell+1}}\bigg)_{\ell=2x-1}\;,\cr}
$$
and for $x\geq 2$, we find an upper bound
$$
0<R_9(x)<{1\over 4x}\sum_{\ell=1}^{+\infty}{1\over 2\ell(\ell+1)}
+{1\over 2}(2x-2){1\over (4x)^2}<{3\over 16x}.
$$
Thanks to an explicit calculation of $U'(x)$ for $x=1{,}2{,}3$,
we get the estimate
$$
|U'(x)|<{0.206\over x}.\leqno(5.9)
$$
Combining (2.13), (3.7), (4.19), (4.22), (4.23), (4.24) and (5.6 -- 5.9), 
we now obtain
$$
\Delta(x)={e^{-4x}\over 4\pi x}\bigg(
-{5\sqrt{2\pi}\over 12\,x^{1/2}}+R(x)\bigg)\leqno(5.10)
$$
with
$$
R(x)=-U'(x)
+\pi\Big(R_1(x)+R_2(x)-R_3(x)\Big)-{5\over 4}R_4(x)+2\,R_5(x)-R_6(x)+R_7(x)
+2\,R_8(x),
$$
whence
$$
|R(x)|<{10.835\over x}.
\leqno(5.11)
$$
These estimates imply (0.10 -- 0.13). The proof of the Theorem is complete.\qed
\bigskip

{\bigbf References}
{
\parindent=14.5mm
\medskip

\item{\st[Art31]} {\petcap E.\ Artin} -- {\it The Gamma Function}~;
Holt, Rinehart and Winston, 1964, traduit de~: {\it Einf\"uhrung in die 
Theorie der Gammafunktion}, Teubner, 1931.

\item{\st[BJ15]} {\petcap R.P.\ Brent \&\ F.\ Johansson} -- {\it  A bound
for the error term in the Brent-McMillan algorithm}~;
Math.\ of Comp., v.~{\bf 84}, 2015, p.~2351--2359.

\item{\st[BM80]} {\petcap R.P.\ Brent \&\ E.M.\ McMillan} -- {\it Some new 
algorithms for high-precision computation of Euler's constant}~;
Math.\ of Comp., v.~{\bf 34}, January 1980, p.~305--312.

\item{\st[Dem85]} {\petcap J.-P.\ Demailly} -- {\it Sur le calcul num\'erique 
de la constante d'Euler}~; Gazette des Math\'ematiciens, vol.\ {\bf 27}, 
(1985), p.~113--126.

\item{\st[Eul14]} {\petcap L.\ Euler} -- {\it De progressionibus harmonicis
observationes}~; Euleri Opera Omnia, Ser.\ 1, v.~{\bf 14},
Teubner, Leipzig and Berlin, 1925, p.~93--100.

\item{\st[Eul15]} {\petcap L.\ Euler} -- {\it De summis serierum numeros 
Bernoullianos involventium}~; Euleri Opera Omnia, Ser.\ 1, v.15, 
Teubner, Leipzig and Berlin, 1927, p.~91--130. See in particular
p.~115.\ The calculations are detailed on p.~569--583.

\item{\st[Swe63]} {\petcap D.W.\ Sweeney} -- {\it On the computation of Euler's 
constant}~; Math, of Comp., v.~{\bf 17} (1963), p.~170--178.

\item{\st[Wat44]} {\petcap G.N.\ Watson} -- {\it A Treatise on the Theory of
Bessel Functions}~; 2${}^{\rm nd}$  edition, Cambridge Univ.\ Press,
London, 1944.
\bigskip\bigskip}

Jean-Pierre Demailly\\
Universit\'e de Grenoble Alpes, Institut Fourier\\
UMR 5582 du CNRS, 100 rue des Maths\\
38610 Gi\`eres, France
\bigskip\bigskip

\centerline{(October 2016, revised in November 2017)}

\bye